\documentclass[12pt]{article}
\usepackage{amsfonts}
\begin{document}
\def\R{\mathbb{R}}
\def\C{\mathbb{C}}
\def\Z{\mathbb{Z}}
\def\N{\mathbb{N}}
\def\Q{\mathbb{Q}}
\def\D{\mathbb{D}}
\def\Ex{\mathbb{E}}
\def\P{\mathbb{P}}
\def\T{\mathbb{T}}
\def\Sp{\mathbb{S}}
\def\hb{\hfil \break}
\def\ni{\noindent}
\def\x{{\xi}^{\ast}}
\def\i{\indent}
\def\a{\alpha}
\def\b{\beta}
\def\e{\epsilon}
\def\d{\delta}
\def\De{\Delta}
\def\g{\gamma}
\def\qq{\qquad}
\def\L{\Lambda}
\def\E{\cal E}
\def\G{\Gamma}
\def\F{\cal F}
\def\Hi{\cal H}
\def\K{\cal K}
\def\O{\cal O}
\def\A{\cal A}
\def\B{\cal B}
\def\L{\cal L}
\def\M{\cal M}
\def\N{\cal N}
\def\Om{\Omega}
\def\om{\omega}
\def\s{\sigma}
\def\t{\theta}
\def\th{\theta}
\def\Th{\Theta}
\def\z{\zeta}
\def\p{\phi}
\def\m{\mu}
\def\n{\nu}
\def\l{\lambda}
\def\Si{\Sigma}
\def\q{\quad}
\def\qq{\qquad}
\def\half{\frac{1}{2}}
\def\hb{\hfil \break}
\def\half{\frac{1}{2}}
\def\pa{\partial}
\def\hb{\hfil \break}
\def\ni{\noindent}
\def\i{\indent}
\def\half{{1 \over 2}}
\def\Om{\Omega}
\def\om{\omega}
\begin{center}
{\bf PREDICTION THEORY IN CONTINUOUS TIME}
\end{center}
\begin{center}
{\bf N. H. BINGHAM}
\end{center}
\begin{center}
To John Kingman
\end{center}
\ni {\bf Abstract} \ We consider prediction theory for stationary stochastic processes in continuous time.  We discuss prediction using the whole (infinite) past, and using only a finite section of the past.  The solutions to both these classical problems have long been known.  Our focus is to provide short simple proofs which reveal the probabilistic meaning of the results. \\ 

\ni {\it Keywords} \ Cram\'er representation, Kolmogorov isomorphism theorem, stationary process, orthogonal increments, moving-average representation, Szeg\H{o}'s condition, Gaussian process, Paley-Wiener theorem \\
\ni 2010 Mathematics Subject Classification: Primary 60-02 Secondary 62-02\\  

\ni {\bf 1. Introduction} \\

\i We begin with a (weak-sense) stationary (complex-valued) stochastic process $X = \{ X_t: t \in \R \}$, briefly, `stationary process'.  We recall the two principal prediction problems for $X$: for $t > 0$, to find the best (least-squares) prediction $\hat X_t$ given (a) the whole past, $\{ X_u: u \leq 0 \}$, (b) a finite section of the past, $\{ X_u: -2T \leq u \leq 0 \}$ for $T > 0$.  \\
\i The solution to (a) is due to Wiener in 1947 [Wie] (see \S 5 for details  of this and other references, and comments).  It also appears in Doob's classic of 1953, [Doo].  The solution to (b) is due to Krein in 1954 [Kre].    Its solution in the Gaussian case is one of the main results of the 1976 monograph by Dym and McKean [DymM3].    \\
\i The solution of (a) in [Doo] is admirably clear.  Unfortunately, it seems to have been overlooked in the extensive subsequent literature.  As both the authors of [DymM3], and Marcus in his review of it [Mar], comment, the solution of (b) presented there lacks probabilistic meaning.  It is also long and difficult.  Our main object here is two-fold: to present Doob's short, simple proof of (a), extended here to (b), and to make the probabilistic meaning transparent.  \\
\i This note complements our studies of the one-, finite- and infinite-dimensional cases in discrete time [Bin2,3,4], and its method is close to that used there. \\
\i It emerges that to give a short and probabilistically revealing treatment of this material, there are four essential ingredients.  Below, we label these acronymically: \\
(i) the {\it Cram\'er Representation}, $(CR)$, \\
(ii) the {\it Kolmogorov Isomorphism Theorem}, $(KIT)$, \\
(iii) the {moving-average representation}, $(MA)$, \\
(iv) the {\it Szeg\H{o} condition}, $(Sz)$. \\
The second and fourth of these (both distributional) are necessary for any treatment of this area, and are accordingly found in all the sources we cite.  The first and third (both pathwise) are present in Doob's book (see [Cra1,2,3] and the comments in [Doo, 637]).  They seem not to have been used systematically together in the context of prediction theory since [Doo], and that only for problem (a).  Unfortunately this has had a detrimental effect on the development of the area.  Our aim here is to redress this.  See \S 5 for comments on the rather unusual history here. \\      

\ni {\bf 2. Prediction from the whole past} \\

\i We follow Doob [Doo, XII], to whose notation and terminology we largely conform, for ease of reference. \\
{\it 1.  Cram\'er Representation and Kolmogorov Isomorphism Theorem} \\
\i We take a {\it stationary} process $X$ as above.  In continuous time, some regularity condition is needed on the paths to avoid pathology (recall {\it Belyaev's dichotomy} [Bel]: for a Gaussian process, the paths are, a.s., either continuous or pathological -- unbounded above and below on every interval).  We take them continuous in mean square.  We recall the {\it Cram\'er Representation} (or {\it spectral representation})
$$
X_t = \int e^{2 \pi it\m} dY(\m)                             \eqno(CR)
$$
(integrals are over $\R$ unless otherwise stated), where the process $Y = \{ Y(\m): \m \in \R \}$ has {\it orthogonal increments}.  If the covariance function is $r$, then by Khinchin's result of 1934 [Khi], $r$ is positive definite, and so a Fourier transform: 
$$
r(t) = \int e^{2 \pi it\m} dG(\m),
$$
where $G$, the {\it spectral measure}, is a complex-valued Lebesgue-Stieltjes measure function of finite variation.  Then
$$
{\Ex} [ \vert dY(\m) \vert^2] = dG(\m),                        \eqno(Spec)
$$
and the {\it Kolmogorov Isomorphism Theorem} holds:
$$
X_t \leftrightarrow e^{it.} := [\m \mapsto e^{it\m}].           \eqno(KIT)
$$
See e.g. [Doo, IX-XII],  [CraL, \S 7.5], [Lin, \S 3.3.2]. \\
{\it 2.  Regularity and determinism; Szeg\H{o}'s condition} \\
\i The (best-possible) prediction error ${\s}^2(t)$ of $X(t)$ at time $t > 0$ given the whole past $\{ X(u): u \leq 0 \}$ may be zero.  In this case, which is in some sense degenerate (and unrealistic: there is `no room for new randomness'), the whole process is completely determined by its remote past (`at $- \infty$').  Such a process is called {\it deterministic} (though it is {\it random} as the remote past is!), otherwise it is called {\it regular}.  The condition for regularity is {\it Szeg\H{o}'s condition}
$$
\int \frac{\log G'(\m)}{1 + {\m}^2} d\m > - \infty              \eqno(Sz)
$$
[Doo, XII, Th. 5.1]; immediately before this, Doob shows how to obtain this condition from the perhaps more familiar form in discrete time, for which see e.g. [Doo, XII, Th. 4.3], or [Bin2]. \\
\i It emerges that the process $X$ splits into a regular component, for which the prediction involves only the {\it spectral density} (the absolutely continuous component $G'$ of the spectral measure $G$), and a deterministic one, involving the singular component of $G$.  For the second, prediction is not needed, and it passes through the prediction process unchanged.  So to simplify the exposition we restrict attention here to the regular component, assume $X$ is regular, and refer to [Doo] for details of the general case involving both components (see also \S 4.4). \\
{\it 3.   The Szeg\H{o} function and the Fourier transform} \\
\i Recall (see e.g. [PalW], or [Doo, 435]) that for functions $f \in L_2(\R)$, the Fourier transform $f \leftrightarrow f^*$ is `symmetrical': $f^* \in L_2$, and \\
$$
f(t) = \int e^{2 \pi i ts} f^{*}(s) ds, \qquad
f^*(s) = \int e^{-2 \pi i ts} f(t) dt.
$$
\i Under the Szeg\H{o} condition $(Sz)$, Doob [Doo XII, Th. 5.2] shows that there exists a function $c^*(t) \in L_2$, vanishing on $[0.\infty)$ and with (complex) Fourier transform $\int c^*(t) e^{2\pi i wt} dt$ non-zero in the open lower half-plane $Im \ w < 0$, whose (real) Fourier transform $c(\m)$ satisfies
$$
{\vert c(\m) \vert}^2 = G'(\m)
$$
( so $c$ `is the complex square root of the spectral density').  One has
$$
\log \Bigl[ \int c^*(t) e^{2\pi t} dt]\Bigr] 
= \log \Bigl[ \int \frac{c(\m)}{1 - i \m} \Bigr] d\m
= \frac{1}{2 \pi} \Bigl[ \frac{ \log G(\m)}{1 + {\m}^2} \Bigr] d\m
$$
(so $c^*$ decreases exponentially at infinity).  For convenience, we preserve symmetry by calling either or both of $c, c^*$ the {\it Szeg\H{o} function}. \\
{\it 4. Prediction and Hilbert-space geometry} \\
\i In regression, one obtains a best linear unbiased estimator (BLUE) by projection, orthogonality and use of Pythagoras's theorem in Euclidean space.  In this infinite-dimensional setting, one likewise obtains a prediction in the least-squares sense, by projection and use of orthogonality and Pythagoras's theorem in Hilbert space.  See e.g. [Doo, 635-7] (Comments on Ch. X) for historical comments here, in particular on the Russian school where this originated. \\
\i Given a set of random variables, one is thus concerned with their closed linear span $\M \{ . \}$, the subspace they generate.  Here it is convenient to use Doob's notation for the time- and frequency-domain versions [Doo, 581]:
$$
{\N}(r,s):={\M}\{ X(t):t\in (r,s] \},\q
{\N}_t := {\N}(-\infty, t) = {\M}\{X(s): s \leq t \},  
$$
$$
(r,s){\N} := {\M}\{e^{2\pi it\m}:t\in (r,s]\}, \ \
_t{\N} := (-\infty,t) {\N} = {\M}\{e^{2\pi is\m}: s \leq t \}.  
$$
\i In $(KIT)$, random variables correspond to functions; as these will have argument $\m$ here, it is convenient to write this as 
$$
X(t) \leftrightarrow e^{2 \pi it\m},
$$
and more generally, if $\psi$ corresponds to $\Psi(\m)$, 
$$
\psi = \int \Psi(\m) dY(\m).
$$
\ni {\it 5.  Orthogonal-increments processes and the moving-average representation} \\
\i Doob [Doo, XII, Th. 5.2] gives a representation for the process $X$ in terms of the Szeg\H{o} function $c^*$ and an orthogonal-increments process $\xi$, 
$$
X(t) = \int_{-\infty}^0 c^*(s) d\xi (t + s)
     = \int_{-\infty}^t c^*(u - t) d\xi (u).                  \eqno(MA)
$$
We use the name $(MA)$ for {\it moving average}: this is the continuous analogue of a moving-average result in discrete time.  Doob [Doo, Th. X.8] shows that a stationary process in discrete time has a moving-average representation with mutually orthogonal random variables if and only if its spectral measure is absolutely continuous (showing once again the crucial importance here of the spectral density). \\
\i So under $(Sz)$, one now has
$$
X(t) = \int c^*(s) d\xi (t + s), 
= \int e^{2 \pi it \mu} c(\m) d\x (\m),                       \eqno(MA^*)
$$
by Parseval's formula; here the processes $\xi, \x$ are Fourier transforms of each other, have orthogonal increments, and
$$ 
\Ex [ \vert d\xi (t) \vert^2 ] = dt,    \qquad
\Ex [ \vert d\x (\m) \vert^2 ] = d\m,
$$
by Plancherel's theorem.  In the first one has the {\it time domain}, in the second the {\it frequency domain}; the link between the two is $(KIT)$.  Referring to $(CR)$, 
$$
dY(\m) = c(\m) d\x (\m): \qquad d\x (\m) = dY(\m)/c(\m).
$$
{\it 6. Doob's formula for the prediction} \\
\i Doob [Doo, XII, Th. 5.2] further gives: \\
(a) the prediction given the whole past at lag $\tau$ (in which $\int_{-\infty}^t$ in $(MA)$ is replaced by $\int_{-\infty}^{t - \tau}$), which we write as 
$$
\hat{{\Ex}} [X(t) \vert \{ X(u): u \leq t - \tau \}]
= \int_{-\infty}^{-\tau} c^*(s) d\xi (t + s)
= \int_{-\infty}^{t - \tau} c^*(u - t) d\xi (u);         \eqno(Pred)
$$ 
(b) the prediction error, 
$$
{\s}^2(t) = \int_{-\tau}^0 \vert c^*(s) \vert^2 ds.      \eqno(Err)
$$
For, the right-hand side of $(Pred)$ is in ${\N}_{t - \tau}$, while the (actual) prediction error 
$$
X(t) - \hat{{\Ex}} [X(t) \vert \{ X(u): u \leq t - \tau \}]
= \int_{t - \tau}^t c^*(u - t) d\xi(u)
= \int_{-\tau}^0 c^*(s) d\xi(t + s)
$$ 
is orthogonal to  ${\N}_{t - \tau}$ (as it involves only $\xi$-increments with arguments $\geq t - \tau$), showing that $(Pred)$ gives the required orthogonal projection.  \\
\i For (b), the (mean-square) prediction error ${\sigma}^2(\tau)$ is the expectation of the squared modulus of the right-hand side above.  By the It\^o isometry and $\Ex [ \vert d\xi (t) \vert^2 ] = dt$, this is $\int_{-\tau}^0 \vert c^*(s) \vert^2 ds$, as required.  \\
\i The same result is given in `$\m$-function' or`prediction function' language in [Doo] (the first display on the next page, p.590), as a quotient
$$
{\Psi}_t(\m) = e^{2\pi it\m} \int_{-\infty}^{-t} e^{2 \pi i \m s} c^*(s) ds/c(\m).
$$
For as $c(\m) = \int e^{2\pi i \m s} c^*(s) ds$,
$$
e^{2 \pi its} - {\Psi}_t 
= e^{2 \pi it\m} \int_t^{\infty} e^{2 \pi i \m s} c^*(s) ds.
$$
The terms on the left correspond to $X(t)$ and its prediction, which is in $_0{\N}$; their difference, the right-hand side, is orthogonal to $_0{\N}$ as $t > 0$. \\
{\it 7. Wiener's formula for the Wiener filter} \\
\i In the special case when $X$ is {\it Gaussian}, when uncorrelatedness and independence are the same, this result is due to Wiener [Wie].  For then, this last quotient is (to within notation) {\it Wiener's formula} for the {\it Wiener filter} [Wie, (2.0393), (2.041)].  The same formula is also obtained in [DymM3, 3, 88, 90], and in operator language in [AroD, Ch. 9], in higher dimensions.  See \S 5 for comments here. \\      
   
\ni {\bf 3. Prediction from part of the past} \\

\i We turn now to the second of the two classical prediction problems: prediction of $X(t)$ given a finite section of the past at lag $\tau > 0$ of length $2T$, $\{ X(u): u \in [t - \tau -2T, t - \tau] \}$ for $T > 0$.  The problem was briefly mentioned in [Doo, XII.5].  It was solved by Krein in 1954 [Kre], just `post-Doob'.  In the Gaussian case, it is one of the main themes of [DymM3]; see [DymM3, \S 6.10]. \\
\i The method of Doob above is immediately applicable here.  It gives:: \\
(a) the prediction: $(Pred)$ above becomes
$$
\hat{\Ex} [X(t) \vert X(u), u \in [t - \tau -2T, t - \tau] \}
= \int_{t - \tau - 2T}^{t - \tau} c^*(u - t) d\xi(u).  
$$
The actual prediction error is thus
$$
\Bigl( \int_{-\infty}^{t - \tau - 2T} + \int_{t - \tau}^t \Bigr) c^*(u - t) d\xi(u),
$$
as this and the prediction sum to $X(t)$; \\ 
(b) the mean-square error, ${\sigma}^2(\tau, T)$.  For this, take the expectation of the squared modulus of the above.  Of the four resulting integrals, the cross terms vanish by orthogonal increments.  The two squared terms are evaluated by the It\^o isometry as above, giving 
$$
{\sigma}^2(\tau, T) = 
\Bigl( \int_{-\infty}^{-2T - t}+\int_{-t}^0 \Bigr) \ \vert c^*(s) \vert^2 ds. 
$$
{\it Note}. 1.  No special new machinery is in fact needed here. \\
2. Entire functions of exponential type at most $T$, as in the solution in [DymM3], correspond under the Fourier transform to the length $2T$ on which we are given data, by the Paley-Wiener theorem ([PalW, Th. X]; [Boa, \S 6.8], [Koo, Vol. I, IIIC]).\\  
  
\ni {\bf 4. Complements} \\

\ni {\it 1.  The Gaussian case} \\
\i When $X$ is Gaussian, the orthogonal-increments process $\xi$ in $(MA)$ can be taken as Brownian motion.  Equivalently, one can use the language of white noise.  See e.g. [DymM3, \S 4.5]. \\
\ni {\it 2. Hardy spaces}\\
\i  Hardy spaces are extensively used in [DymM3] (and e.g. [Bin2]), and in [Doo] though not by that name. \\
\i The Lebesgue decomposition of the spectral measure and the corresponding decomposition of the Cram\'er Representation hold in general, by $(CR)$ and $(Spec)$.  Under the Szeg\H{o} condition $(Sz)$, the Szeg\H{o} function exists and is an {\it outer} function in Beurling's terminology, and corresponds to the regular component of the process and the absolutely continuous component in the Lebesgue decomposition of the spectral measure.  Inner functions correspond to the deterministic component of the process (and are needed to analyse their structure), and the singular component of the measure.  Doob does not address this, so does not need to use Hardy-space language; Szeg\H{o}'s work is cited as the {\it Szeg\H{o} alternative} in [DymM3, \S 4.2]. \\
\i The decomposition into regular and deterministic components is the {\it Wold decomposition}.  Accordingly, the link above is called the {\it Wold-Cram\'er concordance}.  It is exact in one dimension.  In higher dimensions, it is exact only when the relevant matrix has full rank; see e.g. [Bin3]. \\
{\it 3. Higher dimensions} \\
\i In discrete time, we refer to [Bin1,2,3]; the continuous case can be handled analogously.  The (one- and) finite-dimensional case in continuous time is the main theme of [AroD].  As in [Bin3], the computational side of functional data analysis (FDA) reduces to the finite-dimensional case; the dimension $p$ needs to be chosen appropriately.  We do not pursue the infinite-dimensional continuous-time case (FDA with a continuum of curves) here for reasons of space. \\
\ni {\it 4. Passage between discrete and continuous time} \\
\i  Several methods have been used here.  One is the M\"obius map $w = (z - i)/(z + i)$, which maps the half-plane conformally onto the disc [Bin2, \S 8.11].  A second is replacing discrete time $n$ by continuous time $t$ in the Cram\'er Representation.  The process thus defined is {\it harmonizable} in Lo\`eve's sense [Lo\`e]; see [BinM, \S 2] for properties of the resulting process, links with the sampling theorem, etc.  A third is Doob's elegant reduction of continuous to discrete time in [Doo, 582-4]. \\
{\it 5.  The deterministic component} \\
\i  Though degenerate in some sense, the deterministic component is of both mathematical and practical interest.  First, it is difficult in practice to disentangle it from the regular component; see e.g. [DymM3, 89].  Second, rates of convergence in the deterministic case are of interest; see [Ros], [BabGT]. \\
\ni {\it 6. Germ fields and splitting fields} \\
\i  There are many interesting problems, more detailed and demanding than those considered here, concerning (in particular) Gaussian processes in continuous time.  Germ fields -- on `the interplay between the immediate past and the immediate future', and how they can be split by conditioning, as with the Markov property -- are one example.  We cite here Levinson and McKean in 1964 [LevM], Pitt in 1972 [Pit], and Dym and McKean in 1970 [DymM1,2] and 1976 [DymM3, \S 4.10].  Another is splitting fields -- the smallest subfield of the past conditioning on which makes the rest of the past and the future independent; see [DymM, \S 4.8, 4.9] for the relevant {\it splitting test}, and formula for the {\it splitting type}, the smallest $T$ for which $[-2T,0]$ will serve here. \\
{\it 7. Levinson-McKean and Nehari aspects} \\
\i  Although the work in [LevM] seems quintessentially continuous-time, it is worth noting that ideas in it can be usefully taken over to discrete time [KasB].   \\
\ni {\it 8. The logarithmic integral} \\
\i  Koosis's book [Koo] gives a wealth of applications of the Szeg\H{o} (logarithmic) integral in analysis. \\
\ni {\it 9. Inverse spectral and spectral problems} \\
\i  Problems of this sort are treated at length in [DymM3] (one dimension) and [AroD] (several).  For links between Szeg\H{o} theory and spectral theory, see Simon's books [Sim1,2,3]. \\
\ni {\it 10. Entropy} \\
\i  The Szeg\H{o} condition is clearly of {\it entropy} type.  For background here, see e.g. [Bin2, 301-2] (Gibbs Variational Principle), and more recently, Bessonov and Denisov [BesD1,2], Bessonov [Bes]. \\
\ni {\it 11.  De Branges spaces} \\
\i  During the 1960s, de Branges [deBra] developed his theory of Hilbert spaces of entire functions, now called {\it de Branges spaces}.  For their applications in probability theory, see e.g. [Dym], [DymM3], [AroD]. \\  

\ni {\bf 5. Historical comments} \\

\i The historical development of this area has been unusual.  With the benefit of hindsight, it seems odd that the above was not done in the late 1970s, and that some of the main strands here developed independently and without reference to (or apparent awareness of) each other or the literature.  These comments are offered here as the author's attempt to address this. \\
\ni {\it 1. Wiener} \\ 
\i Wiener's famous book [Wie] -- the `Yellow Peril', written in 1942 during WWII, originally classified, and published in 1949 -- gives, among much else, Wiener's formula for the Wiener filter [Wie, (2.0393), (2.041)].  It mentions Kolmogorov -- see [Wie, 59] for a full account of the connections between Kolmogorov's 1941 paper and his work -- but does not mention Szeg\H{o} or Cram\'er by name. \\
\ni {\it 2. Doob} \\
\i Doob's 1953 classic treats both aspects of the mathematics relevant here, distributional [Doo, XII.5, Th. 5.2] and pathwise [Doo, XII.5, Th. 5.3], as well as giving (not by that name) Wiener's formula for the Wiener filter [Doo, end of XII.5].  He also treats the crucial {\it Cram\'er Representation} (again, not by that name) in [Doo, IX] (processes with orthogonal increments), [Doo, XI] (stationary processes -- continuous parameter) and [Doo, XII] (linear least-squares prediction -- stationary (wide-sense) processes).  We note that the text of [Doo] should be read in conjunction with the Appendix (comments) at the end (it might have been better to give these chapter-wise or even section-wise). \\
\i The penultimate paragraph of the Preface to [Doo] begins `Chapter XII, on prediction theory, is somewhat out of place in the book, since it discusses a rather specialized problem'.  (The author is very glad it is there, and regards it, with the foundational work on separability, measurability and versions (Ch. II), and the chapter on martingales (Ch. VII), as one of the `three crowning glories' of the book.)  The paragraph concludes `I had the benefit of stimulating conversations with Norbert Wiener on this subject.'  More is true: the book had been intended to be a collaboration between them (comment to the author, 1976, quoted in [Bin1, \S 1]).  We mention this here, partly for its historical interest, partly because it may account for (at least part of) the subsequent `parallel lines' development of the area. \\
\i It may also be relevant here that these matters are treated right at the end of a long and demanding book.  For comments on (the 50th anniversary of) the book, see [Bin1]. \\
\ni {\it 3. Krein} \\
\i While prediction given the whole past (back to $-\infty$) stems largely from the work mentioned above of Kolmogorov, Wiener and Cram\'er in the 1940s, prediction given only part of the past (a time-interval $[-2T,0]$ say) is more recent, and stems from work of Krein in 1954 (and so is `post-Doob').  For an assessment of Krein's work on prediction theory, see [Dym]. \\
{\it 4. Dym and McKean} \\
\i The more recent (1976) modern classic of Dym and McKean [DymM3] gives Wiener's formula for the Wiener filter [DymM3, 3, 88, 90, 91].  It cites Szeg\H{o} but not Cram\'er. \\
\i Marcus [Mar] gives a fine and scholarly review of [DymM3].  But, as he remarks in his penultimate paragraph, `However, this is not a book for browsing'.  Just as `you have to mean it to read Doob', `you have to mean it to read Dym and McKean'.  Taken together, these two no doubt account for at least some of the `parallel-lines' development. \\
\i Dym and McKean [DymM3, \S 6.10] give (very briefly) their solution to the problem of prediction from a finite section of the past, illustrating it with a number of examples.  They also solve the problem of {\it interpolation} [DymM3, \S 4.13, 6.13, 6.14]: `filling in' an interval between the past and the future (think of gaps in the fossil, geological or archaeological record). \\
\i The preferred technical tool in [DymM3] is {\it Krein's theory of strings}, stemming from 1954 [Kre].  This, with its emphasis on second-order differential operators of `$d^2/dxdm$' type, is thematically reminiscent of the 1965 classic [It\^oM] by It\^o and McKean.  We note that  reducing the time-interval on which data are given corresponds to `shortening the string' [DymM3, Ch. 6]. \\  
\i Perhaps it is a pity that Cram\'er did not give a textbook account of these matters; see [CraL, 144]. \\
{\it 5. Arov and Dym} \\
\i Arov and Dym, in their recent book of 2018 {\sl Multivariate prediction, de Branges spaces, and related extension and inverse problems} [AroD], inevitably to some extent a sequel to and updating of [DymM3], address prediction (witness their title; their Ch. 9 has the same title, `Past and future', as Ch. 4 of [DymM3]), and the higher- (though finite-)dimensional case.  They cite Doob and Wiener, though, again, not Cram\'er.  As in \S 4.10, the preferred technical tool in [AroD] is de Branges spaces. \\
\ni {\it 6. The author} \\
\i The author resolved on reading [DymM3] in 1976 and [Mar] in 1978 to unearth the probabilistic meaning here.  This required more familiarity with Doob (whence [Bin1]), prediction theory in discrete time (whence [Bin2,3,4]), and Gaussian processes (whence [BinS1,2]).  Meanwhile, as always, there was plenty else to do, hence the long delay. \\
\ni {\it 7. Conclusion} \\
\i This account may seem unduly long, but is included here as the author's best attempt to account for how something so important, useful, and probabilistically interesting could have remained `hidden in plain sight' for 45 years. \\

\ni {\bf Acknowledgements} \\

\i I thank Harry Dym, Henry McKean and Mike Marcus for drawing attention in [DymM3] and [Mar] to the need to give probabilistic meaning to the material there, which led to this paper.  It is a pleasure to recall my debt to Joe Doob, professionally and personally.  It is also a pleasure to quote again Kingman's Dictum [Bin1, \S 3]: `It's all in Doob', and to dedicate this paper to John Kingman, to whom its contents will be no surprise. \\

\ni {\bf References} \\

\ni [AroD] Arov, D. Z. and Dym, H., {\sl Multivariate prediction, de Branges spaces, and related extension and inverse problems}.  Operator Th. Anal. Appl. {\bf 266}, Birkh\"auser, 2018.  \\
\ni [BabGT] Babayan, N. M., Ginovyan, M. S. and Taqqu, M. S., Extensions of Rosenblatt's results on the asymptotic behaviour of the prediction error for deterministic stationary sequences.  {\sl J. Time Series Analysis} {\bf 42} (2021) (Murray Rosenblatt Memorial Issue), 622-652, arXiv:2006.0043.\\
\ni [Bel] Belyaev (Belayev), Yu. K., Continuity and H\"older's conditions for sample functions of stationary Gaussian processes.  {\sl Proc. Fourth Berkeley Symp. Math. Stat. Probab. Volume II: Contributions to Probability Theory} (ed. J. Neyman), 23-33, U. California Press, 1961. \\
\ni [Bes] R. V. Bessonov, Entropy function and orthogonal polynomials.  {\sl J. Approximation Thepry} {\bf 272} (2021), 105650, arXiv:2104.11196.\\
\ni [BesD1] Bessonov, R. V. and Denisov, S. A., A spectral Szeg\H{o} theorem on the real line.  {\sl Adv. Math.} {\bf 359} (2020), 41 p., arXiv:1711.05671. \\
\ni [BesD2] Bessonov, R. V. and Denisov, S. A., Zero sets, entropy, and pointwise asymptotics of orthogonal polynomials.  {\sl J. Functional Analysis} {\bf 280} (2021), 109002, arXiv:1911.11280. \\
\ni [Bin1] Bingham, N. H., Doob, fifty years on. {\sl J. Appl. Probab.} {\bf 42} (2005), 257-266.  \\
\ni [Bin2] Bingham, N. H., Szeg\H{o}'s theorem and its probabilistic descendants.  {\sl Probability Surveys} {\bf 9} (2012), 287-324. \\
\ni [Bin3] Bingham, N. H., Multivariate prediction and matrix Szeg\H{o} theory.  {\sl Probability Surveys} {\bf 9} (2012), 325-339. \\
\ni [Bin4] Bingham, N. H., Prediction theory for stationary functional time series. \\
\ni [BinM] Bingham, N. H. and Missaoui, Badr, Aspects of prediction. {\sl J. Appl. Prob.} {\bf 51A} (2014), 189-201. \\
\ni [BinS1] Bingham, N. H. and Symons, Tasmin L., Gaussian random fields on sphere and sphere cross line.  {\sl Stochastic Proc. Appl.} (Larry Shepp Memorial Issue); arXiv:1812.02103; {\tt https://doi.org/10.1016/j.spa.2019.08.007}. \\
\ni [BinS2] Bingham, N. H. and Symons, Tasmin L.,  Gaussian random fields: with and without covariances.  {\sl Theory of Probability and Mathematical Statistics} (Special Issue in honour of M. I. Yadrenko, ed. A. Olenko), to appear.  \\
\ni [Boa] Boas, R. P., {\sl Entire functions}.  Academic Press, 1954. \\
\ni [Cra1] H. Cram\'er, On the theory of stationary random processes.  {\sl Ann. Math.} {\bf 41} (1940), 215-230 (reprinted in {\sl Collected Works of Harald Cram\'er, Volume II}, 925-940, Springer, 1994). \\
\ni [Cra2] H. Cram\'er, On harmonic analysis in certain function spaces.  {\sl Ark. Mat. Astr. Fys.} {\bf 28B} (1942), 1-7 ({\sl Works II}, 941-947). \\
\ni [Cra3] H. Cram\'er, A contribution to the theory of stochastic processes.  {\sl Proc. Second Berkeley Symposium Mat. Stat. Prob} (ed. J. Neyman) 329-339, U. California Press, 1951 ({\sl Works II}, 992-1002). \\
\ni [CraL] Cram\'er, H. and Leadbetter, M. R., {\sl Stationary and related stochastic processes: Sample function properties and their applications}.  Wiley, 1967. \\
\ni [deBra] de Branges, L., {\sl Hilbert spaces of entire functions}.  Prentice-Hall, 1968. \\
\ni [Doo]  Doob, J. L., {\sl Stochastic processes}.  Wiley, 1953. \\
\ni [Dym] Dym, H., Krein's contributions to prediction theory.  {\sl Operator Theory: Advances and Applications} {\bf 118} (2000), 1-15, Birkh\"auser. \\
\ni [DymM1] Dym, H. and McKean, H. P., Applications of de Branges spaces of integral functions to the prediction of stationary Gaussian processes.  {\sl Illinois J. Math.} {\bf 14} (1970), 299-343. \\
\ni [DymM2] Dym, H. and McKean, H. P., Extrapolation and interpolation of stationary Gaussian processes.  {\sl Ann. Math. Statist.} {\bf 41} (1970), 1817-1844. \\
\ni [DymM3] Dym, H. and McKean, H. P., {\sl Gaussian processes, function theory and the inverse spectral problem}.  Academic Press, 1976. \\
\ni [I\^toM] It\^o, K. and McKean, H. P., {\sl Diffusion processes and their sample paths}.  Grundl. math. Wiss. {\bf 114}, Springer, 1965. \\
\ni [KasB] Kasahara, Yukio and Bingham, N. H., Verblunsky coefficients and Nehari sequences.  {\sl Trans. Amer. Math. Soc.} {\bf 366} (2014), 1363-1378. \\
\ni [Khi] Khinchin, A. Ya., Korrelationstheorie der stationa\"aren stochastischen Prozesse.  {\sl Math. Annalen} {\bf 109} (1934), 604-615. \\
\ni [Koo] Koosis, P., {\sl The logarithmic integral}, I, 2nd ed., Cambridge Univ. Press, 1998 (1st ed. 1988), II, Cambridge Univ. Press, 1992.\\
\ni [Kre] Krein, M. G., On a fundamental approximation problem in the theory of extrapolation and filtration of stationary random processes (Russian).  {\sl Doklady Akad. Nauk SSSR} {\bf 94} (1954), 13-16; English transl., {\sl Amer. Math. Soc. Selected Transl. Math. Stat. Prob.} {\bf 4} (1964), 127-131. \\
\ni [LevM] Levinson, N. and McKean, H. P., Weighted trigonometrical approximation on $R^1$ with application to the germ field of a stationary Gaussian noise.  {\sl Acta Math.} {\bf 112} (1964), 99-143. \\
\ni [L\'evy] P. L\'evy, {\sl Processus stochastiques et mouvement brownien}.  Gauthier-Villars, Paris, 1948 (2nd ed. 1965). \\
\ni [Lin] Lindgren, G., {\sl Stationary stochastic processes, Theory and applications}.  CRC Press, 2013. \\
\ni [Lo\`e] M. Lo\`eve, Fonctions al\'eatoires du second ordre.  Note, p.228-352, in [L\'evy]. \\
\ni [Mar] Marcus, M. B., Review of [DymM].  {\sl Annals of Probability} {\bf 5} (1977), 1047-1051. \\
\ni [PalW] Paley, R. E. A. C. and Wiener, N., {\sl Fourier transforms in the complex comain}.  AMS Colloq. Publ. {\bf XIX}, Amer. Math. Soc., 1934. \\
\ni [Pit] Pitt, L. D., On problems of trigonometric approximation from the theory of stationary Gaussian processes.  {\sl J. Multivariate Anal.} {\bf 2} (1972), 145-161. \\
\ni [Ros] Rosenblatt, M., Some purely deterministic processes.  {\sl Research Division, College of Engineering, New York University}, 1956, 21p. (MR0080397, J. L. Doob). \\
\ni [Sim1] Simon, B., {\sl Orthogonal polynomials on the unit circle.  Part 1: Classical theory}. AMS Colloquium Publications {\bf 54.1}, American Math. Soc., Providence RI, 2005. \\
\ni [Sim2] Simon, B., {\sl Orthogonal polynomials on the unit circle.  Part 2: Spectral theory}.  AMS Colloquium Publications {\bf 54.2}, American Math. Soc., Providence RI, 2005. \\
\ni [Sim3] Simon, B., {\sl Szeg\H{o}'s theorem and its descendants: Spectral theory for} $L^2$ {\sl perturbations of orthogonal polynomials}.  Princeton University Press, 2011.\\  
\ni [Wie] Wiener, N., {\sl Extrapolation, interpolation and smoothing of stationary time series, with engineering applications}.  MIT Press/Wiley, 1949; Martino, 2013. \\

\ni Mathematics Department, Imperial College, London SW7 2AZ, UK; \\ n.bingham@ic.ac.uk \\ 

\end{document}